# СРАВНЕНИЕ ДВУХ СТАТИСТИЧЕСКИХ АЛГОРИТМОВ РЕКОНСТРУКЦИИ ИЗОБРАЖЕНИЙ ПРИ КОЛИЧЕСТВЕННОЙ ОЦЕНКЕ ПАТОЛОГИЧЕСКИХ ОЧАГОВ МЕТОДОМ ГАММА-ЭМИССИОННОЙ ТОМОГРАФИИ

*А. В. Нестерова[1,a], Н. В. Денисова[2,b]*

[1] Институт математики СО РАН им. С.Л. Соболева, Россия, 630090, г. Новосибирск, пр. Коптюга, д. 4
[2] Институт теоретической и прикладной механики СО РАН им. С.А. Христиановича, Россия, 630090, г. Новосибирск, ул. Институтская, д. 4/1
E-mail: [a] a.nesterova@alumni.nsu.ru, [b] nvdenisova2011@mail.ru

Проведен сравнительный анализ двух статистических подходов к решению задачи реконструкции изображений при обследовании пациентов методом однофотонной эмиссионной компьютерной томографии (ОФЭКТ). Выполнено сравнение алгоритма *Ordered Subset Expectation Maximization (OSEM)*, установленного на большинстве современных ОФЭКТ системах, и алгоритма нового поколения Maximum a Posteriori с заданием априорной информации на основе функционала энтропии (MAP-Ent) при количественной оценке накопления радиофармпрепарата (накопленной активности) в патологических очагах. Исследования выполнены методом имитационного компьютерного моделирования с использованием цифрового двойника стандартизированного вещественного фантома NEMA IEC, который включает 6 сфер разного размера, имитирующих очаги поражений. Оценка точности реконструкции патологических очагов проводилась с помощью максимального коэффициента восстановления (Recovery Coefficient – $RC_{max}$), равного отношению значения реконструированной накопленной активности в максимуме к ее истинной величине. Результаты исследований показали, что поведение алгоритма OSEM в итерационном процессе характеризуется неустойчивостью, появлением «шума» и краевых артефактов на изображениях очагов. Использование постфильтрации стабилизирует решение и ведет к его сходимости, что позволяет получить кривую зависимости $RC_{max}$ от размера очага, которая потенциально может применяться в качестве поправки к клиническим оценкам очагов поражений. Однако при этом существенно занижается оценка активности в очагах малого размера и возможна потеря малых очагов на изображениях. Алгоритм MAP-Ent демонстрирует принципиально иное поведение: он обеспечивает устойчивую сходимость и количественную точность оценок активности без необходимости постфильтрации, сохраняя контрастность даже малых очагов. Однако номер итерации, на которой достигается точное решение в каждом очаге, существенно зависит от выбора единственного (глобального) параметра регуляризации, что ограничивает возможность оптимальной реконструкции очагов разного размера. Результаты данных исследований показали необходимость перехода к локальной регуляризации в алгоритме MAP-Ent, позволяющей адаптивно подбирать параметры для обеспечения количественной точности реконструкции очагов поражений.

Ключевые слова: однофотонная эмиссионная компьютерная томография, алгоритмы реконструкции изображений, OSEM, MAP-Ent, регуляризация, количественная оценка.

# COMPARISON OF TWO STATISTICAL IMAGE RECONSTRUCTION ALGORITHMS FOR QUANTITATIVE ASSESSMENT OF PATHOLOGICAL LESIONS USING GAMMA EMISSION TOMOGRAPHY

*A.V. Nesterova[1,a], N. V. Denisova[2,b]*

[1] Sobolev Institute of Mathematics, SB RAS, 4 Koptyug ave., Novosibirsk, 630090, Russia
[2] Khristianovich Institute of Theoretical and Applied Mechanics SB RAS, 4/1 Institutskaya st., Novosibirsk, 630090, Russia

This study compares two statistical approaches to image reconstruction in single-photon emission computed tomography (SPECT). We evaluated the widely used Ordered Subset Expectation Maximization (OSEM) algorithm and the newer Maximum a Posteriori approach with Entropy prior (MAP-Ent) approach in the context of quantifying radiopharmaceutical uptake in pathological lesions. Numerical experiments were performed using a digital twin of the standardized NEMA IEC phantom, which contains six spheres of varying diameters to simulate lesions. Quantitative accuracy was assessed using the maximum recovery coefficient ($RC_{max}$), defined as the ratio of the reconstructed maximum activity to the true value. The study shows that OSEM exhibits unstable convergence during iterations, leading to noise and edge artifacts in lesion images. Post-filtering stabilizes the reconstruction and ensures convergence, producing $RC_{max}$–size curves that could be used as correction factors in clinical evaluations. However, this approach significantly underestimates uptake in small lesions and may even lead to the complete loss of small lesions on reconstructed images. In contrast, MAP-Ent demonstrates fundamentally different behavior: it achieves stable convergence and preserves quantitative accuracy without post-filtering, while maintaining the contrast of even the smallest lesions. However, the iteration number at which accurate reconstruction is achieved depends strongly on the choice of a single global regularization parameter, which limits optimal performance across lesions of different sizes. These results demonstrate the need for locally adaptive regularization in MAP-Ent to improve quantitative accuracy in lesion reconstruction.

Keywords: single-photon emission computed tomography, image reconstruction algorithms, OSEM, MAP-Ent, regularization, quantitative assessment.

# 1. Введение

Однофотонная эмиссионная компьютерная томография (ОФЭКТ) является диагностическим методом ядерной медицины. При обследовании пациенту вводят радионуклидный фармацевтический препарат (РФП), который обычно состоит из двух составляющих: несущей (таргетной) молекулы и радионуклида. Несущая молекула подбирается таким образом, чтобы РФП накапливался в патологических очагах. При этом радионуклид служит «меткой», испускающей гамма-кванты. Гамма-излучение регистрируется из разных положений вращающейся вокруг пациента гамма-камеры. На основе измеренных данных осуществляется решение обратной задачи реконструкции пространственного распределения РФП в организме – источника испускания гамма-квантов.

Полученные ОФЭКТ-изображения существенно отличаются от изображений магнитно-резонансной томографии (МРТ) и рентгеновской компьютерной томографии (КТ), отражающих морфологические свойства различных анатомических структур пациента. ОФЭКТ визуализирует распределение РФП, которое соответствует нормальным или патологическим метаболическим процессам, протекающим в организме человека на клеточном и молекулярном уровнях. ОФЭКТ имеет относительно низкое пространственное разрешение по сравнению с МРТ и КТ изображениями, однако основным преимуществом этого метода является высокая контрастность изображений патологических очагов, что делает метод незаменимым для выявления онкологических заболеваний.

Диагностическая точность метода ОФЭКТ существенно зависит от подхода к решению задачи реконструкции пространственного распределения источников излучения по измеренным данным. Радиационный распад радионуклида описывается пуассоновским процессом, соответственно, зарегистрированные данные подчиняются пуассоновскому распределению. С математической точки зрения, проблема реконструкции ОФЭКТ-изображений относится к классу обратных и некорректных задач с пуассоновскими данными и должна решаться с применением статистической регуляризации [1,2].

Первым статистическим подходом к реконструкции ОФЭКТ-изображений стал вариационный метод максимизации функции правдоподобия для пуассоновского распределения, реализованный в алгоритме *Maximum Likelihood Expectation Maximization* (*MLEM*) [3]. В настоящее время большинство функционирующих в мире ОФЭКТ-установок оснащены его ускоренной версией – итерационным алгоритмом реконструкции *Ordered Subset Expectation Maximization* (*OSEM*) [4]. Однако подход на основе максимизации функции правдоподобия не является регуляризованным, поэтому поведение алгоритма *OSEM* характеризуется неустойчивостью: с ростом числа итераций решение сильно «зашумляется» и страдает от появления краевых артефактов на высококонтрастных границах очагов поражений [5]. Чтобы исключить эти неблагоприятные эффекты применяют итерационную регуляризацию – остановку решения после определенного числа итераций, когда ожидается достижение оптимального решения. Однако критерий достижения оптимального решения не определен, поэтому номер останова итерационного процесса выбирают на основе рекомендаций производителя оборудования или практического опыта. Это приводит к тому, что в различных медицинских центрах используют изображения, полученные при разном числе итераций, что значительно снижает воспроизводимость и возможность анализа многоцентровых исследований. Кроме того, прерывание итерационного процесса не снимает полностью влияние «шума» и краевых артефактов. Дополнительным регуляризирующим фактором является пост-фильтрация реконструированного решения. В результате применения пост-фильтрации осуществляется сглаживание решения и уменьшение влияния «шума» и краевых артефактов, однако при этом значительно занижается количественное решение в очагах поражений малого размера [6].

Для уменьшения ошибок, вызванных влиянием шума, необходимо использовать регуляризированные алгоритмы. В литературе высказывалось мнение, что краевые артефакты можно смягчить с помощью применения регуляризированных алгоритмов реконструкции на основе байесовского подхода *Maximum a Posteriori* (*MAP*) с априорной информацией, способствующей сглаживанию решения [7]. В данной работе выполнены исследования алгоритма *MAP* с заданием априорной информации на основе функционала энтропии (*MAP-Ent*). Проведен сравнительный анализ алгоритмов *OSEM* и *MAP-Ent* в решении задачи получения точной количественной оценки накопленной активности в патологических очагах.

Решение проблемы получения точных количественных оценок накопленной активности в очагах поражений является чрезвычайно актуальной задачей на современном уровне развития ядерной медицины, поскольку эти оценки определяют персонализированную поглощенную дозу и эффективность радионуклидной терапии. Используемые на клинических ОФЭКТ системах алгоритмы *OSEM* не позволяют получать точные решения, поэтому в руководстве Европейской Ассоциации по Ядерной Медицине (*EANM*) [8] рекомендуется использовать эмпирические поправки к полученному решению. Эти поправки называются «коэффициент восстановления (*Recovery Coefficient, RC*)» и определяются отношением реконструированного значения активности к истинной активности РФП в очагах. В клинической практике значения RC получают с помощью экспериментальных измерений с применением стандартизированных вещественных фантомов *NEMA IEC* (*National Electrical Manufacturers Association International Electrotechnical Commission*), включающих сферические вставки, имитирующие очаги разного диаметра [9-11]. Затем полученные величины RC используются в качестве поправки к реконструированным значениям активности в очагах. Однако проведение широких клинических исследований с подобными фантомами ограничены из-за лучевой нагрузки на исследователей и высокой стоимости радиофармпрепаратов. Современной альтернативой является математическое имитационное моделирование и численные эксперименты с использованием цифровых двойников фантомов и сканеров. В данной работе этот метод применяется для сравнительной оценки двух статистических алгоритмов *OSEM* и *MAP-Ent* в решении проблемы получения точной количественной реконструкции активности в очагах поражений.

Настоящая работа является логическим продолжением предыдущего исследования [6], в котором была продемонстрирована ограниченность стандартного алгоритма *OSEM* в достижении точной количественной оценки патологических очагов. Основной вывод работы [6] заключался в необходимости разработки и анализа новых методов реконструкции на основе байесовского подхода *MAP* с применением статистической регуляризации.

## 2. Постановка задачи и методы моделирования

В работе проведено компьютерное имитационное моделирование клинической процедуры ОФЭКТ/КТ исследования вещественного стандартизированного фантома *NEMA IEC* для оценки точности реконструированных изображений патологических очагов и определения поправочных коэффициентов (*Recovery Coefficient, RC*). Целью работы является сравнительный анализ двух алгоритмов: *OSEM* и *Maximum A Posteriori* с энтропийной регуляризацией (*MAP-Ent*) – в точности количественной оценки накопленной активности в очагах поражений. Моделирование проводилось в приближении к реальным клиническим условиям: учитывались параметры сканирования, физические процессы регистрации гамма-квантов и влияние алгоритмов реконструкции на количественные характеристики изображений. В качестве моделируемого радионуклида использовался технеций-99m ($^{99m}$Tc), наиболее широко применяемый в клинической практике изотоп с энергией излучения 140 кэВ.

### 2.1. Физико-математическая модель

Физические основы модели:
1. Радиоактивный распад описывается законом:
$$n(t) = n(0)e^{-\frac{t}{\tau}},$$
где $n(t) = \{n_j(t): j = 1, \dots, J\}$ – распределение концентрации РФП в трехмерном объекте в момент времени $t$; $n(0)$ – начальное распределение РФП; $\tau$ – среднее время жизни радиоактивного изотопа. Для $^{99m}$Tc среднее время жизни составляет около 6 часов, тогда как продолжительность стандартного сбора проекционных данных методом ОФЭКТ $\Delta t$ обычно составляет около 20-30 минут.

Поскольку испускание гамма-квантов носит случайный характер и концентрация РФП низкая, число испущенных гамма-квантов во всех направлениях в единицу времени описывается случайным полем $f = \{f_j: j = 1, \dots, J\}$. При условии $\Delta t \ll \tau$ среднее значение гамма-квантов, испускаемых из единичного объема, считается постоянным за время наблюдения и пропорциональным локальной активности РФП: $\bar{f} \approx n(0)$.

2. При прохождении в средах различной плотности гамма-квант может пройти до детектора без взаимодействия, либо быть рассеянным или поглощенным. Основными видами взаимодействия являются комптоновское рассеяние и фотоэффект.

3. Для реконструкции *3D*-распределения активности по проекционным данным необходимо знать направление прилетевших на детектор гамма-квантов. В модели ОФЭКТ для этого используется коллиматор. Теоретический идеальный коллиматор обеспечивает точное изображение точечного источника, однако реальные коллиматоры формируют размытое изображение: размер получаемого пятна возрастает с удалением источника от поверхности коллиматора. Этот эффект описывается функцией рассеяния точечного источника (ФРТ).

При численном моделировании дискретное представление исследуемого объекта задается на трёхмерной сетке вокселей. Поток гамма-квантов, проходящий через ткани пациента и фиксируемый детекторами, формирует двумерные проекционные данные для каждого угла обзора $g = \{g_i, i = 1, \dots, I\}$.

Связь между распределением активности РФП и проекционными данными задается уравнением:
$$\sum_j a_{ij} f_j = g_i \text{ или } A\mathbf{f} = \mathbf{g},$$

где $a_{ij}$ – это случайный оператор, который описывает, какая часть гамма-квантов, испускаемых из $j$-го вокселя, обнаруживается $i$-м пикселем детектора. Здесь $f_j$ – это ненаблюдаемые пуассоновские случайные величины с неизвестными средними значениями $\bar{f}_j$, а $g_i$ – наблюдаемые данные, пуассоновские случайные величины с неизвестными средними значениями $\bar{g}_i$. Эти средние значения связаны системой линейных уравнений:

$$\sum_j \bar{a}_{ij} \bar{f}_j = \bar{g}_i,$$

где $\bar{a}_{ij}$ – это вероятность того, что гамма-квант, испущенный $j$-м вокселем, будет зарегистрирован $i$-м пикселем детектора. Вероятности $\bar{a}_{ij}$ образуют системную матрицу, которая считается известной и учитывает эффекты ослабления излучения $P_{ik}^{att}$ и прохождения через систему коллиматор-детектор $P_{kj}^{col-det}$:

$$\bar{a}_{ij} = P_{ik}^{att} \cdot P_{kj}^{col-det}.$$

Формулировка задачи: задача реконструкции ОФЭКТ-изображений сводится к определению средних значений $\bar{f}_j$ по заданным $\bar{a}_{ij}$ и измеренным данным $g_i$. Эта задача относится к классу некорректных обратных задач с пуассоновскими данными.

### 2.2. Методы моделирования

Моделирование выполнялось с использованием программного комплекса «Виртуальная платформа для имитационных испытаний метода ОФЭКТ/КТ» [12-14], включающего последовательные этапы: задание распределения активности РФП в математическом фантоме, генерацию проекционных данных, реконструкцию изображений и оценку качества реконструкции.

Распределение активности РФП моделировалось с использование цифрового двойника фантома *NEMA IEC*. Фотография реального фантома приведена на рис. 1, а. Он содержит шесть сфер различного диаметра (37, 28, 22, 17, 13, 10 мм), имитирующих патологические очаги с повышенным накоплением РФП, а также цилиндрическую легочную вставку диаметром 51 мм, расположенную в центре. Для удобства анализа изображений все центры сфер расположены в одной плоскости.

При клинических измерениях фантом заполнялся водным раствором РФП с разной концентрацией в сферах и основной емкости. Цифровой двойник фантома задавался в декартовой системе координат $\{x, y, z\}$ в области, разделенной на 128×128×92 воксела. Трехмерная модель распределения активности $^{99m}$Tc по воксeлам фантома, называемая картой активности, представлена на рис. 1, б. Поперечное сечение этой карты, проходящее через центр всех сфер, показано на рис. 1, в.

Отношение активности сфера/фон (контраст) было установлено равным 10:1. В клинических исследованиях это отношение может быть варьироваться. В нашей предыдущей работе [6] показано, что значения коэффициента восстановления *RC* изменяются лишь незначительно при отношениях 5:1, 10:1, 20:1, 30:1. Случаи низкой накопленной активности в очагах, например, при отношении 2:1, требуют отдельного исследования.

Активность задавалась в относительных единицах: фон – 10, легочная вставка – 0, очаги – 100. По смоделированным проекциям рассчитывалось общее число зарегистрированных квантов. Затем все значения масштабировались с использованием коэффициента, связывающего общее число насчитанных квантов в модельных исследованиях с измеренными клиническими данными (общее число импульсов). Это позволяло переводить активности $\bar{f}_j$, заданные в относительных единицах, в единицы «импульс/воксель», используемые при реконструкции. В данной работе изображения представлены именно в этих единицах.

Кроме карты активности, для цифрового фантома *NEMA IEC* была создана карта ослабления, описывающая коэффициенты ослабления воды и воздуха для гамма-квантов с энергией 140 кэВ, испускаемых радионуклидом $^{99m}$Tc.

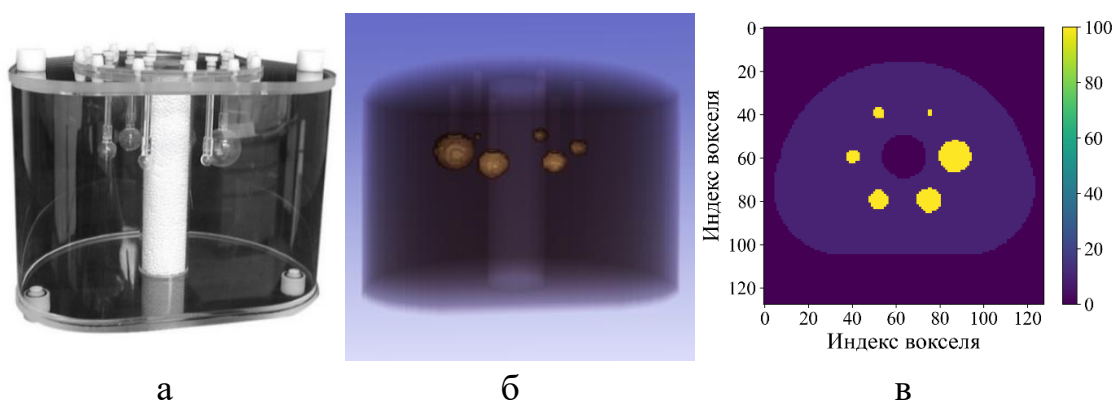

а б в

**Рис. 1.** Сравнение физического фантома *NEMA IEC* и его цифровой модели: (а) фотография физического фантома; (б) трехмерное изображение цифрового двойника фантома; (в) его поперечное сечение, проходящее через центр всех сфер.

Генерация проекционных данных выполнялась с учётом основных физических эффектов, включая ослабление гамма-излучения, а также влияние пространственного разрешения коллиматора. Расчёты проводились для 120 позиций гамма-камеры для круговой орбиты в 360° с шагом 3° и общим числом регистрируемых импульсов ~5 млн.

В рамках программного комплекса «Виртуальная платформа для имитационных испытаний метода ОФЭКТ/КТ» расчет проекционных данных может осуществляться двумя способами: методом Монте-Карло и расчетом дискретизированных уравнений переноса гамма-излучения. Второй подход дает «точные проекции», в которые затем добавляется «пуассоновский шум» с помощью метода исключения Неймана. Оба подхода были верифицированы в предыдущих исследованиях путем сравнения с клиническими данными и тестовыми измерениями.

В данной работе использовался второй подход, поскольку большое число углов сбора данных (120), высокая статистика (~5 млн импульсов) и большое число различных вариантов расчетов потребовали бы огромных компьютерных ресурсов и времени при использовании метода Монте-Карло.

### 2.3. Алгоритмы реконструкции

Исследовались два статистических алгоритма реконструкции изображений *OSEM* и *MAP-Ent*. Математической основой для этих алгоритмов является предположение о пуассоновском распределении проекционных данных. В рамках метода максимального правдоподобия (*MLEM*) решение находится путем максимизации функции правдоподобия, которая для пуассоновской модели имеет вид:

$$P(g|\bar{f}) = \prod_i exp\left(-\sum_j \bar{a}_{ij}\bar{f}_j\right)\frac{\left(\sum_j \bar{a}_{ij}\bar{f}_j\right)^{g_i}}{g_i!}.$$

На практике максимизируют её логарифм, что приводит к итерационной формуле алгоритма *MLEM*:

$$\bar{f}_j^{(n+1)} = \frac{\bar{f}_j^{(n)}}{\sum_i \bar{a}_{ij}} \sum_i \frac{g_i \bar{a}_{ij}}{\sum_k \bar{a}_{ik} \bar{f}_k^{(n)}}.$$

Алгоритм *OSEM* является стандартизированным алгоритмом реконструкции ОФЭКТ и ПЭТ изображений и представляет собой модификацию метода *MLEM*, ускоряющую сходимость за счёт обновления решения по подмножествам проекционных данных (*subsets*) $S_b$:

$$\bar{f}_j^{(n+1)} = \frac{\bar{f}_j^{(n)}}{\sum_{i \in S_b} \bar{a}_{ij}} \sum_{i \in S_b} \frac{g_i \bar{a}_{ij}}{\sum_k \bar{a}_{ik} \bar{f}_k^{(n)}}.$$

Для характеристики вычислительной работы алгоритма *OSEM* и сравнения различных протоколов реконструкции руководство *EANM* [8] рекомендует использовать универсальный параметр: «число обновлений» (*number of updates*), определяемый как произведение количества подмножеств на количество итераций. В данной работе реконструкция проводилась с 10 подмножествами, а число итераций для алгоритма *OSEM* варьировалось в пределах от 1 до 4-х. Соответственно, число обновлений составляло от 10 до 40.

Байесовский подход (*MAP*) состоит в максимизации апостериорной вероятности $P(\bar{f}|g)$, которая, согласно теореме Байеса, пропорциональна произведению функции правдоподобия и априорной вероятности:

$$P(\bar{f}|g) \sim P(g|\bar{f}) \cdot P(\bar{f}).$$

Максимизация логарифма апостериорной вероятности приводит к уравнению:

$$\frac{\partial}{\partial \bar{f}_j}\left[\ln P(g|\bar{f}) + \ln P(\bar{f})\right] = 0.$$

Развитие байесовского подхода на основе принципа энтропии для решения обратных и некорректных задач со стохастическими данными было выполнено в работах [15,16]. В работе [16], функционал энтропии, который используется для задания плотности априорной вероятности, был получен в логарифмическом виде:

$$\ln P(\bar{f}) = \ln P(m) = -\beta \sum_j \left(\bar{f}_j \ln \frac{\bar{f}_j}{m_j} - \bar{f}_j + m_j\right),$$

где коэффициент $\beta$ является параметром регуляризации, $m_j$ описывает априорную информацию до проведения измерений. Разработка и применение этого подхода для реконструкции ОФЭКТ изображений представлены в работе [17]. Подстановка данного априорного распределения и функции правдоподобия для пуассоновского распределения $P(g|\bar{f})$ в уравнение максимизации и его решение приводит к выводу итерационной формулы *MAP-Ent* алгоритма. В качестве априорной информации $m$ на каждом шаге использовалось решение предыдущей итерации $\bar{f}^{(n-1)}$:

$$\bar{f}_j^{(n+1)} = \bar{f}_j^{(n)} exp\left[\gamma \sum_i \left(g_i \frac{\bar{a}_{ij}}{\sum_k \bar{a}_{ik}\bar{f}_k^{(n)}} - \bar{a}_{ij}\right)\right],$$

где параметр $\gamma = \frac{1}{\beta} > 0$ контролирует соотношение между вкладом априорной вероятности (функционал энтропии) и измеренными пуассоновскими данными (функционал правдоподобия).

### 2.4. Поправочный коэффициент к количественной оценке реконструкции активности в очагах поражений

Для количественной оценки ОФЭКТ изображений очагов рассчитывался коэффициент восстановления (*Recovery Coefficient, RC*):

$$RC_{max} = \frac{A_{rec,max}}{A_{true}},$$

где $A_{\text{rec,max}}$ – максимальное значение активности в реконструированном очаге, $A_{\text{true}}$ – истинное значение в очаге. Значения *RC_{max}*, близкие к 1, соответствуют точному восстановлению активности в максимуме. *RC_{max}* меньше 1 означает недооценку, а *RC_{max}* больше 1 – переоценку активности. Следует заметить, что оценка качества реконструкции по одному вокселю (с максимальной интенсивностью) вызывает немало вопросов, однако она используется в клинических исследованиях, поскольку для расчета среднего значения по очагу *RC_{mean}* требуется определить контур очага, что существенно усложняет оценку и может приводить к ошибкам.

### 3. Результаты вычислительных экспериментов и обсуждение
### 3.1. Количественный анализ алгоритма OSEM

Вычислительные эксперименты проводились с применением алгоритма реконструкции *OSEM* как без постфильтрации, так и с постфильтрацией фильтром Баттерворта с параметрами, используемыми в клинической практике. На рис. 2 представлены соответствующие реконструированные поперечные сечения, проходящие через центры сфер, полученные с 1-ой по 4-ю итерацию. Ниже представлены профили активности для сфер 13 и 10 мм (верхний ряд), 17 и 37 мм (средний ряд) и 22 и 28 мм (нижний ряд). Результаты без постфильтрации характеризуются выраженным шумовым искажением и краевыми артефактами. Это связано с отсутствием регуляризации в алгоритме *OSEM*, что при увеличении числа итераций усиливает шумовую составляющую [18]. Краевые артефакты проявляются в виде ложных пиков активности.

Применение постфильтрации стабилизирует решение, однако контрастность изображений малых очагов (13, 17 мм) значительно снижается, что может приводить к неопределенности в диагностике, либо полной потере этих очагов на изображениях.

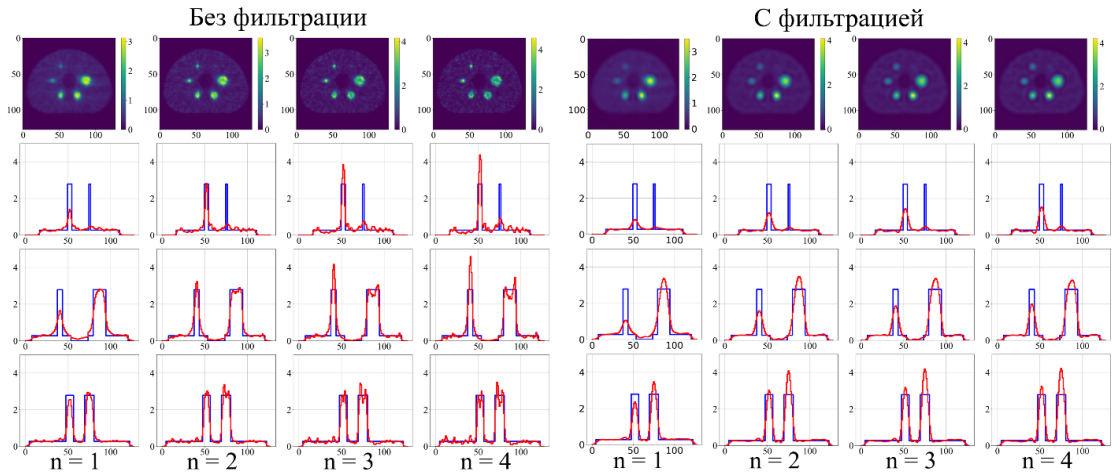

**Рис. 2**. Поперечные сечения (верхний ряд) и соответствующие горизонтальные профили активности (нижние ряды), реконструированные алгоритмом *OSEM*. Сравнение проведено при числе итераций n = 1–4 с постфильтрацией и без нее. Красная кривая – профиль, полученный после реконструкции; синяя – исходный профиль цифровой модели. Профили построены по линиям, проходящим через центр верхних (10 и 13 мм), средних (17 и 37 мм) и нижних (22 и 28 мм) сфер. По осям графиков профилей: горизонтальная ось – индекс вокселя; вертикальная ось – значение активности в условных единицах.

Для изображений, реконструированных алгоритмом *OSEM* с постфильтрацией и без нее, был рассчитан коэффициент восстановления $RC_{max}$. На рис. 3 представлены кривые $RC_{max}$ в зависимости от диаметра сфер для различных значений числа итераций.

Анализ кривых $RC_{max}$, полученных для изображений без фильтрации, показывает, что при увеличении числа итераций разброс между кривыми сохраняется или усиливается. Это свидетельствует о высокой чувствительности решения *OSEM* без постфильтрации к числу итераций и делает невозможным формирование надежной оценки поправочного коэффициента $RC_{max}$ для очагов разных диаметров.

Введение постфильтрации, напротив, приводит решение к сходимости, что позволяет останавливать алгоритм после определенного фиксированного числа итераций. Однако, как было указано выше, корректнее использовать не просто число итераций, а число обновлений. Следует отметить, что применение постфильтрации сопровождается снижением точности реконструкции малых очагов (13 и 17 мм): их профили сильно занижены, а $RC_{max}<1$, что указывает на недооценку интенсивности в очаге поражения.

Было исследовано влияние соседства больших очагов на оценку $RC_{max}$ для малых очагов. В численных экспериментах рассматривались фантомы с тремя очагами диаметром 10, 13 и 17 мм. Полученные значения $RC_{max}$, рассчитанные для изображений с постфильтрацией, практически совпали с аналогичными значениями для фантома с шестью очагами. Таким образом, наличие крупных очагов не оказывает заметного влияния на оценку $RC_{max}$ для малых очагов при установленных расстояниях между сферами в фантоме *NEMA IEC*.

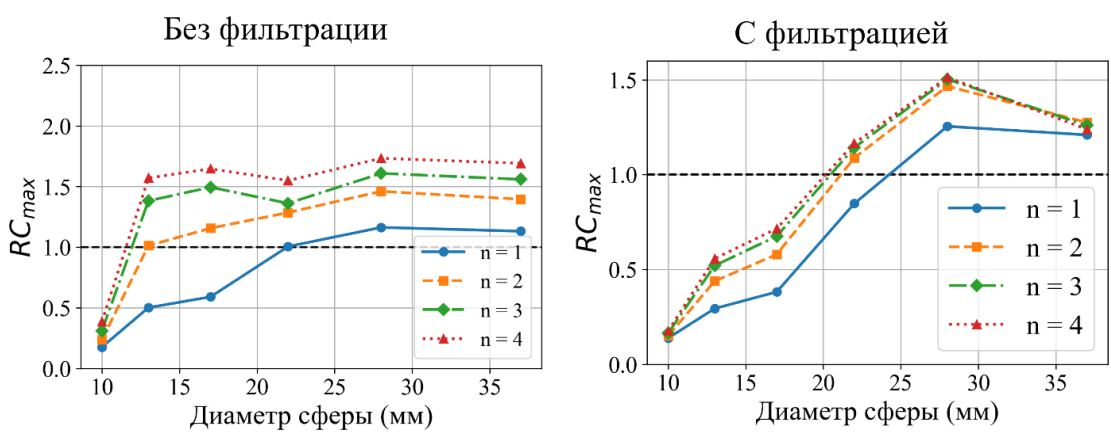

**Рис. 3.** Значения коэффициента восстановления $RC_{max}$ реконструированной активности алгоритмом *OSEM* в зависимости от диаметра сфер. Показаны данные для итераций $n = 1$–4.

Следует отметить, что количественные оценки накопленной активности в очагах на изображениях без фильтрации ненадёжны, поскольку изображения и, соответственно, поправочные коэффициенты $RC_{max}$ подвержены влиянию шума и краевых артефактов. Для практической количественной оценки очаговых поражений необходимо использовать поправочные коэффициенты $RC_{max}$, рассчитанные для изображений, полученных с применением постфильтрации. Однако на реконструированных изображениях с постфильтрацией малые очаги могут вызывать неуверенность или полностью исчезать.

Таким образом, несмотря на то что применение постфильтрации позволяет получать более устойчивые количественные оценки, данный подход сопровождается риском потери информации о малых очагах. В связи с этим остаётся актуальным поиск альтернативных методов, которые обеспечивали бы надёжность количественной оценки без существенного ухудшения визуализации. В литературе одним из перспективных направлений является использование теоретически рассчитанных поправок, независимых от алгоритма реконструкции. К таким подходам относятся разработка универсальных уравнений для прогнозирования *RC* [19] и фундаментальных моделей коррекции эффекта частичного объёма (краевых артефактов), учитывающих форму патологического очага и разрешение системы [20].

### 3.2. Количественный анализ алгоритма MAP-Ent

В алгоритме *MAP-Ent* используется параметр γ. Он умножает слагаемое, соответствующее функции правдоподобия. Таким образом, при больших значениях γ метод приближается к алгоритму *OSEM*, что приводит к усилению шума и краевых артефактов. В настоящем исследовании использовались фиксированные значения γ = 0.01 и γ = 0.1.

На рис. 4 представлены результаты реконструкции для $\gamma = 0.01$ и $\gamma = 0.1$ при различных значениях числа итераций *n*. Параметр $\gamma$ можно интерпретировать как аналог «шага по времени» в итерационном процессе: при одинаковых значениях $\gamma \cdot n$ получаем одинаковые решения.

На рис. 5 представлена зависимость коэффициента восстановления $RC_{max}$ для $\gamma = 0.01$ и $\gamma = 0.1$ при разных значениях числа итераций. Видно, что при $\gamma = 0.1$ значения $RC_{\max} \approx 1$ достигаются практически для всех сфер, кроме самой маленькой, но при разном числе итераций. В клинической практике количество метастатических очагов разных размеров у одного пациента может достигать нескольких десятков и получение точных количественных оценок накопленной активности в каждом очаге, используя изображения после разного числа итераций, является проблематичным.

Более эффективным является подход с использованием разных значений $\gamma$, что соответствует идее локальной регуляризации. Сравнение реконструкций при $\gamma = 0.01$ и $\gamma = 0.1$ (рис. 4 и 5) показывает, что при одинаковых значениях произведения ($\gamma \cdot n$) решения практически совпадают. Это означает, что подбор значений $\gamma$ для очагов разного размера позволяет получать оптимальное решение и коэффициенты восстановления $RC_{\max} \approx 1$ на одной и той же итерации.

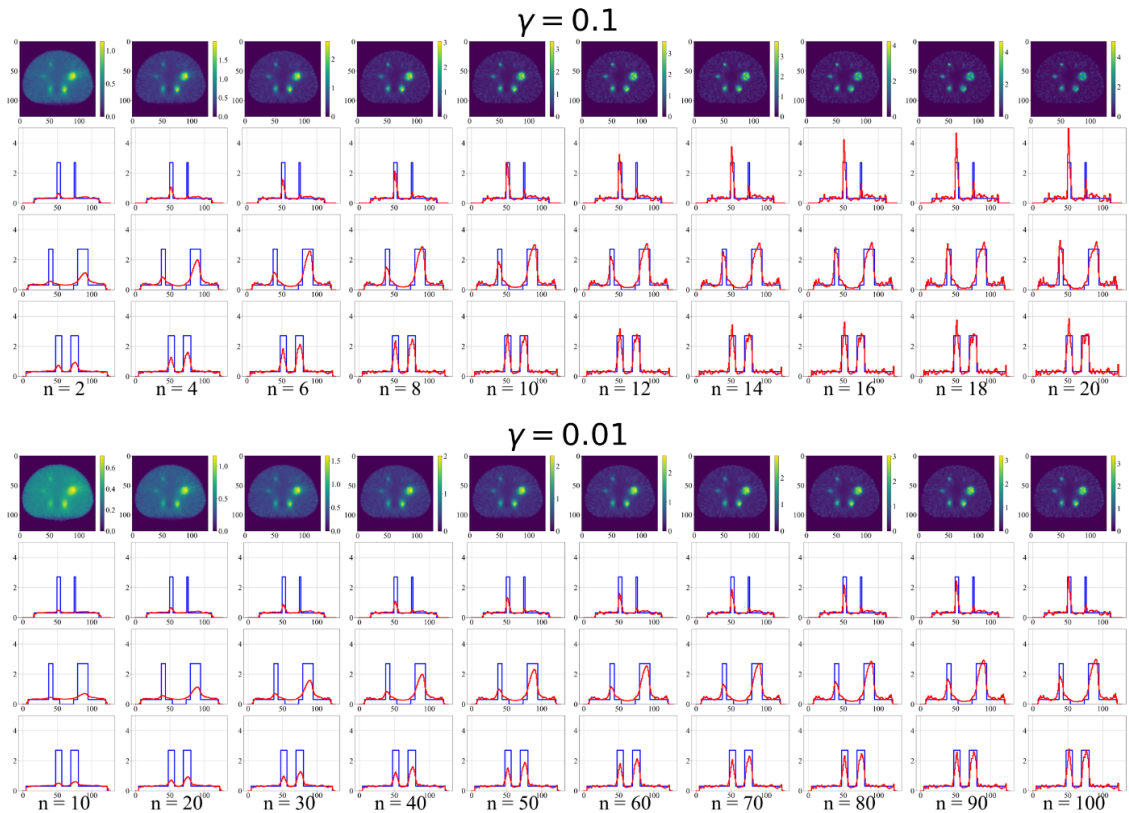

**Рис. 4.** Реконструированные алгоритмом *MAP-Ent* поперечные сечения (верхний ряд) и соответствующие горизонтальные профили активности (нижние ряды). Сравнение проведено для $\gamma = 0.01$ и $\gamma = 0.1$ при числе итераций $n = 1–10$ и $n = 1–100$. Красная кривая – профиль, полученный после реконструкции, синяя – исходный профиль из цифровой модели. Профили построены по линиям, проходящим через центр верхних (10 и 13 мм), средних (17 и 37 мм) и нижних (22 и 28 мм) сфер. По осям графиков профилей: горизонтальная ось – индекс вокселя; вертикальная ось – значение активности в условных единицах.

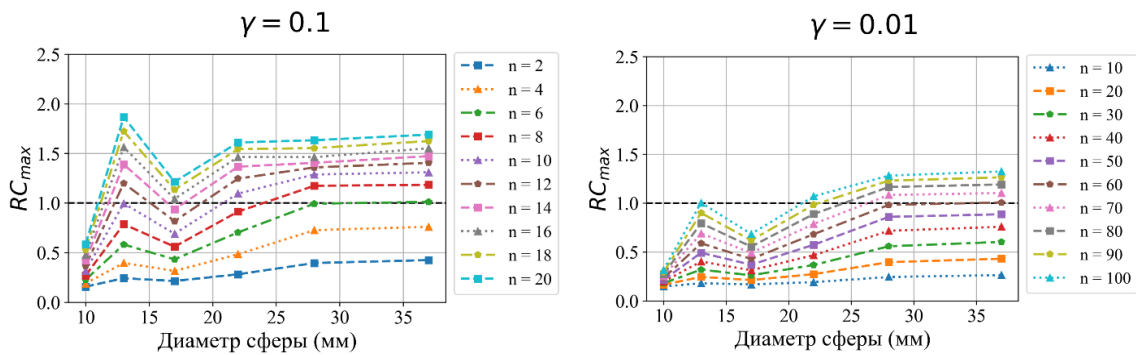

**Рис. 5**. Значения коэффициента восстановления ($RC_{max}$) реконструированной активности алгоритмом *MAP-Ent* для значений параметра $\gamma = 0.01, 0.1$ в зависимости от диаметра сфер. Показаны данные для итераций $n = 1–100$ и $1–20$.

В качестве примера рассчитаем локальные значения $\gamma$ такие, чтобы $RC_{max}$ достигала значения, близкого к 1 на одной и той же итерации для всех очагов. Из рис. 5 видно, что при $\gamma = 0.1$ значение $RC_{\max} \approx 1$ достигается для очагов 37 и 28 мм на 6-й итерации, для сферы 22 – на 9-й, для сферы 17 мм – 15-й итерации. Следовательно, если для сферы 22 мм использовать $\gamma = (9/6) \cdot 0.1 = 0.15$ и для сферы 17 мм использовать $\gamma = (15/6) \cdot 0.1 = 0.25$, то следует ожидать достижения $RC_{\max} \approx 1$ одновременно для всех сфер уже на 6-й итерации. На рис. 6 представлены результаты реконструкции и расчета параметра $RC_{max}$ при использовании локальных значений параметра $\gamma$. Такой подход обеспечивает точную количественную оценку накопленной активности в очагах и потенциально может поднять процедуру ОФЭКТ на более высокий уровень достоверности, поэтому идея применения локальной регуляризации требует дальнейшего развития.

### 3.3. Сравнительная оценка и перспективы развития

Сравнение показывает, что алгоритм *MAP-Ent* превосходит *OSEM* по устойчивости к шуму и количественной точности без необходимости в постфильтрации. Использование *MAP-Ent* позволяет более точно реконструировать распределение активности в очагах без этапа постобработки, обеспечивая количественную достоверность, необходимую для планирования радионуклидной терапии. Алгоритм *OSEM*, в свою очередь, может эффективно применяться для рутинной диагностики, но требует аккуратной настройки параметров и использования постфильтрации для достижения приемлемого качества изображений.

Следует отметить, что для регуляризированных алгоритмов на основе байесовского подхода *Maximum a Posteriori* в литературе отсутствует критерий выбора параметра регуляризации (в нашем случае параметр $\gamma$), и его подбор осуществляется эмпирически. В предыдущей работе одного из авторов [17] был продемонстрирован эффект от использования локальной регуляризации для точной реконструкции очагов поражений с разной интенсивностью накопленной активности. Результаты настоящих исследований продемонстрировали необходимость локальной регуляризации при одинаковой активности, но разных размеров очагов. Такой подход потребует сегментации изображения, что может быть реализовано на основе предварительной реконструкции с использованием *OSEM*. Необходимы дальнейшие исследования этой проблемы.

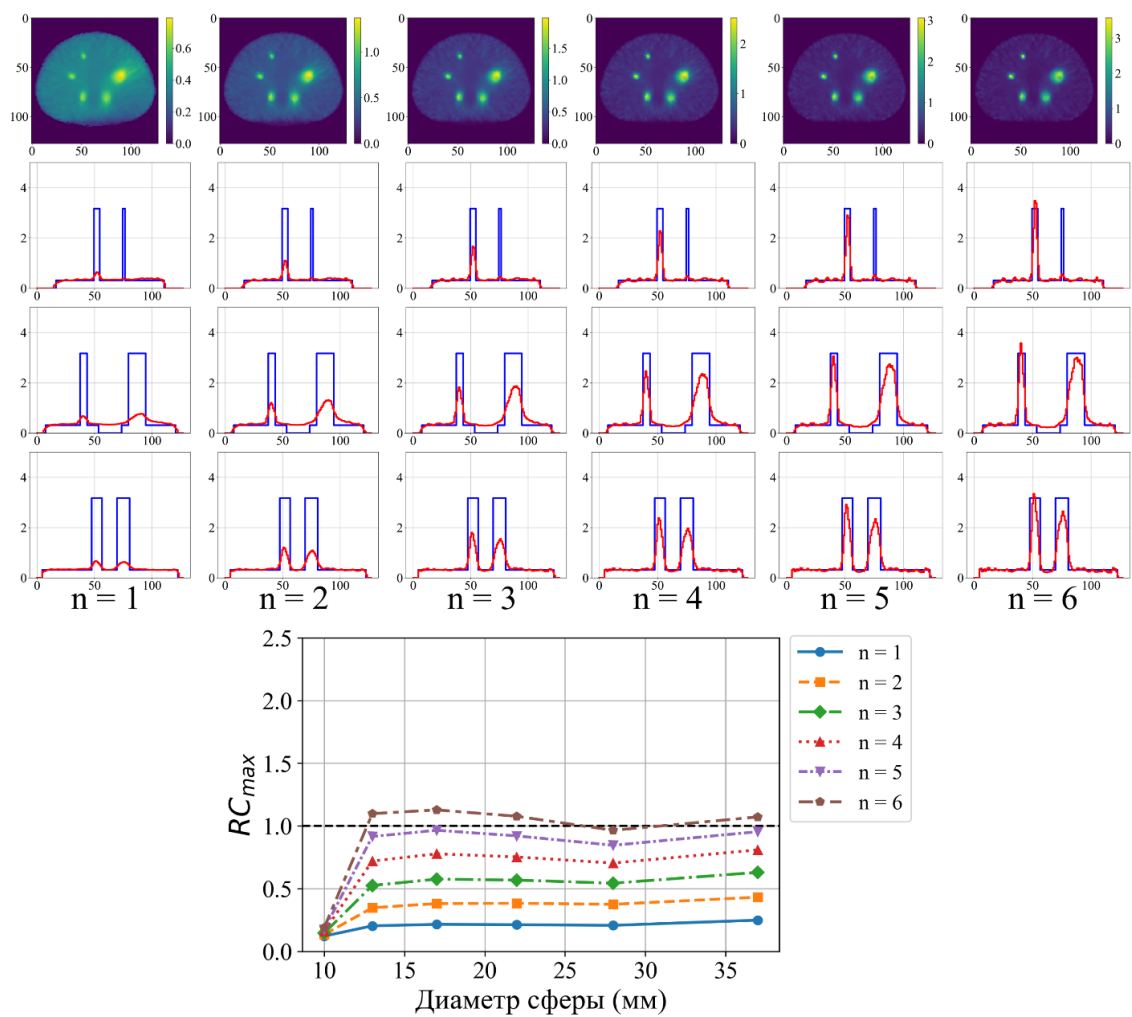

**Рис. 6.** Верхний ряд – реконструированные алгоритмом *MAP-Ent* поперечные сечения (верхний ряд) и соответствующие горизонтальные профили активности (нижние ряды). Сравнение выполнено для различных значений итераций $n = 1–10$ и $n = 1–100$. Красная кривая – профиль, полученный после реконструкции, синяя – исходный профиль из цифровой модели. Профили построены по линиям, проходящим через центр верхних (10 и 13 мм), средних (17 и 37 мм) и нижних (22 и 28 мм) сфер. По осям графиков профилей: горизонтальная ось – индекс вокселя; вертикальная ось – значение активности в условных единицах. Нижний ряд – значения коэффициента восстановления ($RC_{max}$) реконструированной активности алгоритмом *MAP-Ent* в зависимости от диаметра сфер. Показаны данные для итераций $n = 1–6$.

## 4. Заключение

Получение достоверных количественных оценок накопленной активности в очагах поражений является одной из самых актуальных задач современной диагностической ядерной медицины. Для обеспечения точности таких оценок при исследованиях методом ОФЭКТ в алгоритмах реконструкции необходимо учитывать эффекты, влияющие на формирование изображений.

Одним из важных эффектов является размытие изображения точечного источника, что требует учета ФРТ. Это позволяет улучшить визуализацию и повысить возможность обнаружения небольших патологических очагов. Однако на практике учет ФРТ сопровождается появлением краевых артефактов, что приводит к переоценке активности в очагах поражений малого размера. Исследователи-клиницисты пришли к выводу о необходимости постсглаживания изображений для уменьшения таких ошибок реконструкции [21]. В своем комментарии к работе [21] известный специалист в области разработки алгоритмов реконструкции ОФЭКТ и ПЭТ изображений Йохан Нуйтс отметил: «Этот вывод несколько противоречит здравому смыслу, поскольку постсглаживание ухудшает разрешение, в то время как изначально модель учета ФРТ была направлена именно на его улучшение» [7]. Чтобы «смягчить» влияние краевых артефактов в работе [7] было предложено развивать регуляризированные алгоритмы реконструкции на основе байесовского подхода *MAP* с априорными распределениями, сглаживающими решения.

Результаты моделирования, проведенного в настоящей работе, позволили уточнить выводы клиницистов о необходимости постсглаживания (постфильтрации). Действительно, при использовании алгоритма *OSEM* количественные оценки активности без фильтрации оказываются ненадёжными, поскольку изображения и, соответственно, поправочные коэффициенты $RC_{max}$ подвержены влиянию краевых артефактов. Для практического использования необходимо рассчитывать поправочные коэффициенты $RC_{max}$ именно по изображениям с постфильтрацией. С другой стороны, в численных экспериментах получены результаты, показывающие противоречивость применения фильтрации, на что указывал Нуйтс: небольшой очаг диаметром 13 мм, отчетливо различимый на изображении без фильтрации, практически исчез после применения фильтра (рис. 2).

Вторая часть данной работы посвящена исследованию регуляризированного метода *MAP-Ent*. Проведенное исследование показало, что алгоритм *MAP-Ent* обеспечивает устойчивую сходимость и точные количественные оценки активности для очагов различного диаметра. В отличие от стандартного алгоритма *OSEM*, *MAP-Ent* не требует постфильтрации и сохраняет точность даже при высоком контрасте отношения очаг/фон.

Исследования также показали, что поведение решения *MAP-Ent* в итерационном процессе определяется величиной, равной произведению числа итераций и параметра $\gamma$. Это позволяет сделать вывод о возможности использования локального параметра $\gamma$ для достижения точного решения одновременно для всех очагов. Такой подход соответствует методу локальной регуляризации. Применение локальной регуляризации с адаптивным выбором параметра γ для опухолевых очагов представляется перспективным направлением, однако критерии оптимального подбора γ в настоящее время не существуют, что определяет направление дальнейших исследований.

СПИСОК ЛИТЕРАТУРЫ